## ON THE  GENERALIZED RIEMANN  HYPOTHESIS

### Dorin Ghisa

*York University, Glendon College, Toronto, Canada*
*E-mail address:  dghisa@yorku.ca*

### Abstract

It is known that for any primitive character $\chi$ the Dirichlet  $L$-function $L(s\,;\chi)$ verifies a functional equation relating $L(s\,;\chi)$ and $L(1-s\,;\overline{\chi})$ such that $s_0$ is a nontrivial zero of $L(s\,;\chi)$ if and only if $1-\bar{s}_0$ is a nontrivial zero of $L(s\,;\chi)$. Then, proving that this can happen only if $s_0 = 1-\bar{s}_0$, i.e. $\operatorname{Re} s_0 = 1/2$ is equivalent to proving the generalized Riemann Hypothesis for this class of functions.

We show that such an equality is a corollary of some global mapping properties of the functions $L(s\,;\chi)$.

### 1. Introduction

The method we use in this paper for the study of the zeros of the Dirichlet $L$-functions is based on a technique of revealing fundamental domains of these functions. This method has been applied successfully for different classes of analytic functions starting with Blaschke products [ 5 ], [7] - [ 9 ], arbitrary rational functions [6 ], elementary transcendental functions [3] and ending with the Euler Gamma function and the Riemann Zeta function [ 4 ]. This last function is a particular Dirichlet $L$-function and it is understandable that some theorems on arbitrary Dirichlet $L$-functions will duplicate those proved for the Riemann Zeta function. On the other hand, the whole theory will come as a reinforcement of the results we obtained directly when dealing with the Riemann Zeta function and in particular those regarding the location of its zeros [10] - [13].

A *Dirichlet character* $\chi_j$ modulo $q$, $1 \le j \le \varphi(q)$, where $\varphi$ is Euler's totient function, can be defined axiomatically in the following way:

1. $\chi_j(n+q) = \chi_j(n)$, $n \in \mathbb{Z}$

2. If $\gcd(n,q) > 1$, then $\chi_j(n) = 0$, if $\gcd(n,q) = 1$, then $\chi_j(n) \ne 0$

3. $\chi_j(mn) = \chi_j(m)\chi_j(n)$  for all integers $m$ and $n$

By axiom 3, we have $\chi_j(1) = \chi_j(1 \cdot 1) = \chi_j(1)\chi_j(1)$ and since $\gcd(1,q) = 1$ we have by axiom 2 that $\chi_j(1) \ne 0$, hence

4. $\chi_j(1) = 1$ for every modulus $q$ and every $j$, $1 \le j \le \varphi(q)$.

The axiom 1 says that every Dirichlet character $\chi_j$ modulo $q$ is a periodic function defined on $\mathbb{Z}$ with the period $q$. This means that

5. If $a \equiv b(\mathrm{mod}\, q)$, then $\chi_j(a) = \chi_j(b)$

By Euler's Theorem, if $\gcd(a,q) = 1$, then $a^{\varphi(q)} \equiv 1(\mathrm{mod}\, q)$, therefore, by the properties 5 and 4 we have that $\chi_j(a^{\varphi(q)}) = \chi_j(1) = 1$ and by the axiom 3, $\chi_j(a^{\varphi(q)}) = \chi_j(a)^{\varphi(q)}$, so if $\gcd(a,q) = 1$, then $\chi_j(a)^{\varphi(q)} = 1$, and consequently:

6. If $a$ and $q$ are relatively prime numbers, then $\chi_j(a)$ is a root of order $\varphi(q)$ of the unity, for every $j$, $1 \leq j \leq \varphi(q)$.

On the other hand, by axiom 2, if $a$ and $q$ are not relatively prime, then $\chi_j(a) = 0$. In other words, the range of any Dirichlet character is in the set formed with $0$ and the roots of different orders of unity. It can be easily seen that $\overline{\chi_j}$ defined by $\overline{\chi_j}(a) = \overline{\chi_j(a)}$ is a

Dirichlet character if $\chi_j$ is one.

A *Dirichlet series* is a function of the form:

$$(1) \qquad \alpha(s) = \sum_{n=1}^{\infty} a_n n^{-s},$$

where $a_n$ are complex numbers and also $s = \sigma + it$ belongs to some complex domain. It is known (see [14], page 11) that $\alpha(s)$ has an *abscissa of convergence* $\sigma_c \leq +\infty$ such that $\alpha(s)$ converges for all $s = \sigma + it$ with $\sigma > \sigma_c$ and diverges for all $s$ with $\sigma < \sigma_c$. Moreover, if $\alpha(s)$ converges at $s_0 = \sigma_0 + it_0$, then there is a neighborhood of $s_0$ where $\alpha(s)$ converges uniformly. There is a number $\sigma_a$, *the abscissa of absolute convergence*, $\sigma_c \leq \sigma_a \leq \sigma_c + 1$, such that $\alpha(s)$ converges absolutely for all $s = \sigma + it$ with $\sigma > \sigma_a$, but for no $s = \sigma + it$ with $\sigma < \sigma_a$. The strip $\sigma_c \leq \sigma \leq \sigma_a$ of conditional convergence is never wider than $1$ and the value $1$ is reached for the series $\sum_{n=0}^{\infty}(-1)^n n^{-s}$, for which $\sigma_a = 1$ and $\sigma_c = 0$.

A *Dirichlet L-series* is a Dirichlet series

$$(2) \qquad L(s\, ; \chi) = \sum_{n=1}^{\infty} \frac{\chi(n)}{n^s},$$

where $\chi$ is a Dirichlet character $\chi_j$ of some modulus $q$. It is known (see[14], p. 121) that any Dirichlet $L$-series converges for $\mathrm{Re}\, s > 0$. Also, any Dirichlet $L$-series can be extended to an analytic function in the whole plane, except possibly at $s = 1$, where it can have a simple pole. Such a function is called *Dirichlet L-function.* When $q = 1$, the corresponding $L$-function is the Riemann Zeta function. Otherwise, the $L$-series are obtained by replacing $1$ at the numerator of the terms of the Riemann Zeta series by $0$ when $(n,q) \neq 1$ and by some $\varphi(q)$-th roots of the unity, which are the values of the corresponding Dirichlet characters when $(n,q) = 1$.

If $a_n = f(n)$ in the formula (1), where $f$ is a *totally multiplicative function* (see [14], page 20), then an Euler type of product formula:

$$(3) \qquad \alpha(s) = \prod_p \frac{1}{1 - \frac{f(p)}{p^s}},$$

where $p$ runs through the prime numbers, is true in the half-plane of absolute convergence of $\alpha(s)$.

The axiom 3 of the Dirichlet characters shows that for $\operatorname{Re} s > 1$

(4) $\qquad L(s\,;\chi) = \prod_p (1 - \frac{\chi(p)}{p^s})^{-1},$

By using either formula (2) or (4) it can be easily seen that $L(\bar{s}\,;\chi) = \overline{L(s\,;\overline{\chi})}$ for $\operatorname{Re} s > 1$, and this equality holds for the extended functions, hence $L(\bar{s}\,;\chi) = 0$ iff $L(s\,;\overline{\chi}) = 0$. In particular, when $\chi$ is real, then $L(\bar{s}\,;\chi) = 0$ iff $L(s\,;\chi) = 0$, which is a known property when $L(s\,;\chi)$ is the Riemann Zeta function. The product $\chi\overline{\chi}$ is a particular Dirichlet character $\chi_0$, called *principal character and* having all the non zero values equal to $1$.

Multiplying $\zeta(s)$, by $(1 - p^{-s})$, where $p$ are primes dividing $q$, we have for $\operatorname{Re} s > 1$:

(5) $\qquad L(s\,;\chi_0) = \sum_{n=1,(n,q)=1}^{\infty} n^{-s} = \zeta(s)\prod_{p|q}(1 - p^{-s}),$

for any principal character $\chi_0$ modulo $q$. Since the product multiplying $\zeta(s)$ is an integer function, the formula (5) defines $L(s;\chi_0)$ as a meromorphic function in the complex plane having the unique simple pole $s = 1$ with the residue $\prod_{p|q}(1 - \frac{1}{p})$. The zeros of the functions $L(s\,;\chi_0)$ and $\zeta(s)$ are the same, except for the zeros of the product in (5), which are all on the imaginary axis and are considered to be *trivial zeros* of $L(s; \chi_0)$, therefore $\zeta(s)$ and $L(s;\chi_0)$ have the same *non trivial zeros.*

More general, let $q$ be an arbitrary integer, let $d|q$ and denote by $\chi^*$ a character modulo $d$, then $\chi(n) = \chi^*(n)$ if $(n,q) = 1$ and $\chi(n) = 0$ otherwise defines a character modulo $q$ for which

(6) $\qquad L(s;\ \chi) = L(s;\ \chi^*)\prod_{p|q}(1 - \frac{\chi^*(p)}{p^s})$

We say that $\chi^*$ induces $\chi$. A character $\chi^*$ is called *primitive* if it is not induced by any other character.

The formula (5) guarantees (see [ 13 ] ) that a Dirichlet $L$-function corresponding to a Dirichlet principal character verifies the Riemann Hypothesis. Thus, for the purpose of generalizing the Riemann Hypothesis it is enough to deal in what follows only with Dirichlet $L$-functions defined by non principal Dirichlet characters. On the other hand, the formula (6) shows that the zeros of $L(s;\ \chi)$ are the same as those of $L(s;\ \chi^*)$ except for the zeros of the product, which are all on the imaginary axis. They are considered trivial zeros for $L(s;\ \chi)$, hence the non trivial zeros of $L(s;\ \chi)$ and $L^*(s;\ \chi)$ coincide.

The alternative notation $L(q,j,s)$ is used when it is necessary to indicate that $L(s\,;\chi)$ is defined by the $j$-th Dirichlet character $\chi$ modulo $q$. The location of the zeros of two $L$-functions defined by two imprimitive characters $\operatorname{mod} 14$, one real and one complex, is shown in Fig. 1 and 2 below. Taking the pre-image by those $L$-functions of the real axis and coloring differently the pre-image of the negative and that of the positive real half axis, the zeros appear as the junction points of the two colors. The trivial zeros are those situated on the negative real half axis and on the imaginary axis, while the non trivial ones are those situated in the critical strip $0 < \operatorname{Re} s < 1$. Some of them can be seen here in boxes [-5,5] X [-20,20] on the critical line $\operatorname{Re} s = 1/2$. For the purpose of space economy, in these pictures and in some which will follow, the complex plane is rotated counterclockwise by $\pi/2$. It is obviously crucial to study the pre-image of the real axis when dealing with problems related to the location of the zeros of any complex function.

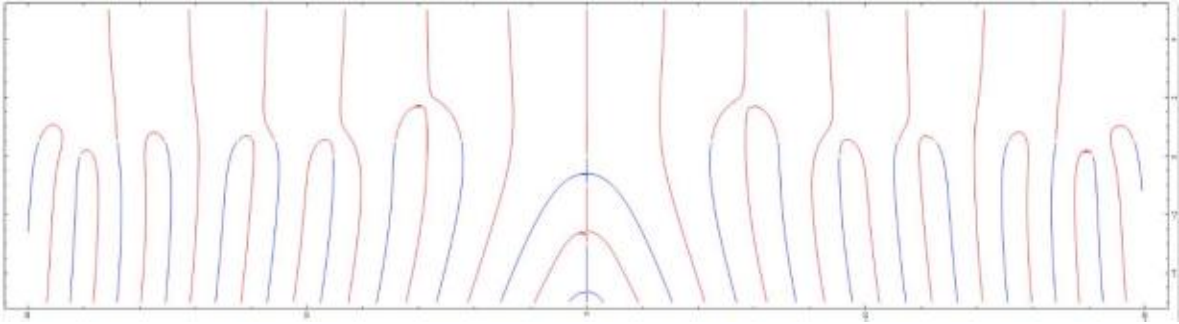

Fig.1    The pre-image of the real axis by $L(14, 4, s)$

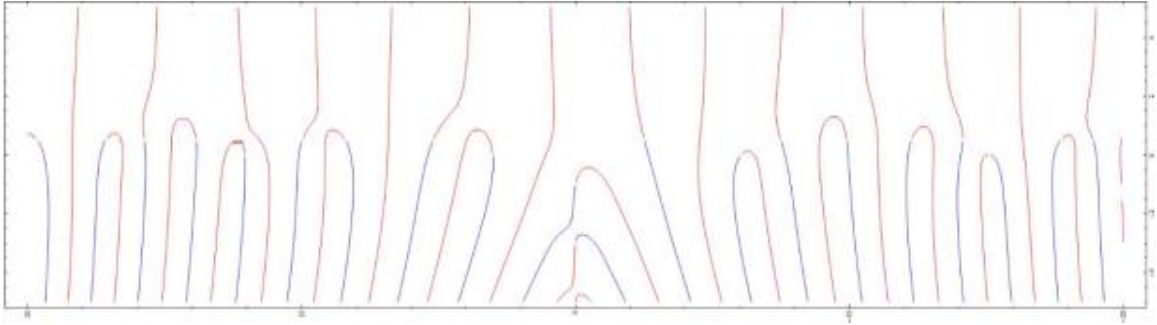

Fig.2    The pre-image of the real axis by $L(14, 2, s)$

The formula (6) suggests that, in what concerns the non trivial zeros, we can deal only with primitive characters. It is known (see [14], Ch.10) that when $\chi$ is a non principal Dirichlet character, then $L(s; \chi)$ is an integer function and if $\chi$ is a primitive character (mod $q$) with $q > 1$, then it verifies a functional equation of the form (see [14], p.333):

(7)     $L(s;\ \chi) = \varepsilon(\chi) L(1 - s;\ \overline{\chi}) 2^s \pi^{s-1} q^{1/2-s} \Gamma(1 - s) \sin \frac{\pi}{2}(s + \kappa),$

where $\varepsilon(\chi) = \frac{\tau(\chi)}{i^\kappa \sqrt{q}}$ with $\kappa = \kappa(\chi) = 0$ if $\chi(-1) = 1$ and $\kappa = 1$ if $\chi(-1) = -1$

The functional equation (7) shows that when $\chi$ is a primitive character, $L(s_0;\ \chi) = 0$ iff $L(1 - \bar{s}_0\ ;\ \chi) = 0$. If $s_0 = \sigma_0 + it_0$, then $1 - \bar{s}_0 = 1 - \sigma_0 + it_0$. Therefore $L(s\ ;\ \chi)$ verifies the Riemann hypothesis iff for every non trivial zero $s_0$, we have $s_0 = 1 - \bar{s}_0$, i.e. $\sigma_0 = 1 - \sigma_0$, or $\sigma_0 = 1/2$. In other words, proving that $L(s\ ;\ \chi)$ verifies the Riemann hypothesis comes to showing that for every two non trivial zeros $s_1 = \sigma_1 + it_1$ and $s_2 = \sigma_2 + it_2$ of $L(s\ ;\ \chi)$, $t_1 = t_2$ implies $\sigma_1 = \sigma_2 = 1/2$. We will show that $t_1 = t_2$ and $\sigma_1 \neq \sigma_2$ represents a contradiction. This fact will result from some global mapping properties of the functions $L(s;\ \chi)$.

## 2. Global Mapping Properties of the Dirichlet $L$-Functions

Finding fundamental domains of the Riemann Zeta function [4], [13] has been an important step in proving the Riemann Hypothesis. The respective domains were obtained by using the pre-image of the real axis. We call *component* of the pre-image by $L(s; \chi)$ of the real axis any maximal curve of that pre-image, in the sense that it has no larger continuation (see [2], p. 23) along the real axis. The components of the pre-image of $\mathbb{R}$ are mapped bijectively by $L(s; \chi)$ onto some intervals of $\mathbb{R}$, or onto $\mathbb{R}$. We will deal also with pre-images of some curves or domains.

The pre-image by $L(s; \chi)$ of an enough small circle $\gamma_r$ centered at the origin and of radius $r$ consists of components which are, in any bounded region, closed curves around every zero of $L(s; \chi)$. As $r$ increases, the respective components expand and $m$ of them ($2 \leq m < \infty$) can touch for an $r = r_0$ at some point $v_0$. In a small neighborhood $V$ of $v_0$ the function $L(s; \chi)$ takes every value on $L(V; \chi) \cap \gamma_{r_0}$ $m$ times, namely once on every component, therefore $v_0$ is a branch point of order $m$ of $L(s; \chi)$. Then $L'(v_0; \chi) = 0$. When increasing $r$ past $r_0$ those components fuse into one closed curve containing the respective zeros. Increasing $r$ more, the fusion process continues and it is possible to obtain unbounded components of the pre-image of $\gamma_r$.

**Theorem 1**. *For any Dirichlet L-function, there is at least one unbounded connected component of the pre-image of the unit disc.*

*Proof:* Let us notice that the formula (2) implies

$$(8) \qquad \lim_{\sigma \to +\infty} L(\sigma + it \; ; \; \chi) = 1$$

for every Dirichlet character $\chi$. It means that for every $\epsilon > 0$ there is $\sigma_t$ such that $\sigma \geq \sigma_t$ implies $|L(\sigma + it \; ; \; \chi) - 1| < \epsilon$. Suppose that all the connected components of the pre-image of the unit disc were bounded. Let us show that in such a case there is $\epsilon > 0$ with the property that all the connected components of the pre-image of the disc centered at the origin and having the radius $1 + \epsilon$ are also bounded. Indeed, for an $r > 0$, let us cover the pre-image $\Gamma$ of the unit circle with the union $O$ of all the open discs centered at points $s \in \Gamma$ and having the radius $r$. The image by $L(s; \chi)$ of these discs are open sets covering the unit circle in the $z$-plane, (where $z = L(s; \chi)$ ), which is a compact set. By the Lebesgue covering theorem a finite covering can be chosen, the pre-image of which still covers $\Gamma$. There is $\epsilon > 0$ such that the circle centered at the origin and of radius $1 + \epsilon$ is included in that finite covering of the unit circle. Therefore the pre-image of this circle is included in $O$, hence all its components are bounded, which contradicts the equality (8) and the assumption we made cannot be true.

Such an unbounded component $D$ of the pre-image of the unit disc contains one or several zeros of $L(s; \chi)$ and for $s \in \partial D$, we can let $\operatorname{Re} s \to +\infty$. The closure of $D$ intersects the critical strip and since the intersection is compact, $D$ can contain only a finite number of non trivial zeros.

**Theorem 2**. *Every unbounded component of the pre-image by $L(s; \chi)$ of the unit disc contains a unique unbounded component of the pre-image of the interval $[0, 1)$ and is included in a strip $S_k$ bounded by two consecutive components $\Gamma'_k$ and $\Gamma'_{k+1}$ of the pre-image of real axis which are mapped bijectively by $L(s; \chi)$ onto the interval $(1, +\infty)$.*

*Proof:* Let us find first the pre-image of a ray $\eta_\alpha$ making a small angle $\alpha$ with the positive real half axis. If $s_0$ is a zero of $L(s; \chi)$, then $L(s_0; \chi) = 0 \in \eta_\alpha$ and if $z = L(s; \chi) \in \eta_\alpha$, then we cannot have $\operatorname{Re} s \to +\infty$, since $z = L(s; \chi)$ keeps away from 1 on $\eta_\alpha$. On the other hand, the continuation by $L(s; \chi)$ along $\eta_\alpha$ starting from $s_0$ is possible, as long as the path over $\eta_\alpha$ does not meet a branch point of $L(s; \chi)$. These are the zeros of $L'(s; \chi)$, hence a discrete set of points whose images can be avoided by $\eta_\alpha$. Therefore we can chose $\alpha$ such that the continuation from any zero $s_0$ of $L(s; \chi)$ along $\eta_\alpha$ is unlimited. Hence $s \to \infty$ on all the components of the pre-image of those $\eta_\alpha$ but $\operatorname{Re} s \nrightarrow +\infty$. For any $\epsilon > 0$, taking $\alpha$ small enough, the arc of the unit circle of opening $2\alpha$ around $z = 1$ is included in a disc of radius $\epsilon$ centered at $z = 1$. Then $\eta_\alpha$ intersects this arc, which means that there is a component of the pre-image of $\eta_\alpha$ intersecting the unbounded component of the pre-image of the unit circle at a point $s'$. The equality (8) implies that $\operatorname{Re} s' \to +\infty$ as $\alpha \to 0$, which shows that indeed, at the limit when $\alpha \to 0$, the respective continuation give rise to two unbounded curves: one contained in the pre-image of the unit disc and which is projected by $L(s; \chi)$ onto the interval $[0, 1)$ and the other denoted $\Gamma'_{k+1}$ exterior to this pre-image and which is projected by $L(s; \chi)$ onto the interval $(1, +\infty)$. The same thing happens on the ray $\eta_{-\alpha}$ giving rise to the same continuation along $[0, 1)$ as $\alpha \to 0$ and to a continuation $\Gamma'_k$ along $(1, +\infty)$. The continuation along $(-\infty, 0)$ is possible as well from any non trivial zero of $L(s; \chi)$ and it is unlimited, since as we'll see later, it cannot meet branch points or poles, so we obtain a unique component $\Gamma_{k,0}$ of the pre-image of the interval $(-\infty, 1)$ which is included in the strip $S_k$ formed by $\Gamma'_k$ and $\Gamma'_{k+1}$. Indeed, $\Gamma_{k,0}$ cannot meet $\Gamma'_k$ or $\Gamma'_{k+1}$ due to the fact that $L(s; \chi)$ is an uniform function.

Fig. 3 below illustrates the creation of a strip $S_k$. The curve $\Gamma_{k,0}$ is that containing a zero $s_{k,0}$ with $\operatorname{Im} s_{k,0} \approx 50$. The strip $S_k$ contains also the curves $\Gamma_{k,-1}$ and $\Gamma_{k,1}$ exhibiting the zeros $s_{k,-1}$ with $\operatorname{Im} s_{k,-1} \approx 48$ and respectively $s_{k,1}$ with $\operatorname{Im} s_{k,1} \approx 53$.

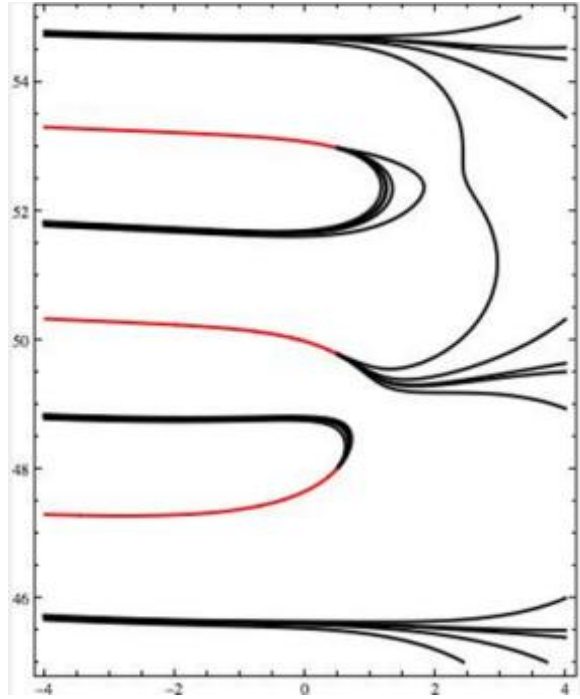

Fig. 3    The pre-image of $\eta_\alpha$ give rise to a strip $S_k$ as $\alpha \to 0$

The following affirmations are corollaries of the well known *argument principle*. If a component $C_k$ of the pre-image of the unit circle contains $j_k$ zeros $s_k$ of $L(s;\ \chi)$, then as $s$ travels on $C_k$ the point $L(s;\ \chi)$ travels on the unit circle turning $j_k$ times around the origin. On its way on $C_k$, the point $s$ meets $j_k$ points $u_{k,j}$ such that $L(u_{k,j};\ \chi) = 1$ when $C_k$ is bounded and $j_{k-1}$ such points when $C_k$ is unbounded and in both cases it meets $j_k$ points $u'_{k,j}$ such that $L(u'_{k,j};\ \chi) = -1$. The points $u_{k,j}$ and $u'_{k,j}$ alternate. Every half-open arc between consecutive points $u_{k,j}$ or consecutive points $u'_{k,j}$ is mapped bijectively by $L(s;\ \chi)$ onto the unit circle. The same thing happens with the unbounded parts of $C_k$, when they exist, starting at the first, respectively ending at the last $u_{k,j}$. The components $\Gamma_{k,j}$, $j \neq 0$ of the pre-image of the real axis are all mapped bijectively by $L(s;\ \chi)$ onto the whole real axis. Every component $\Gamma_{k,j}$, $j \neq 0$ of the pre-image of $\mathbb{R}$ intersects $C_k$ at the points $u_{k,j}$ and $u'_{k,j}$ and the component $\Gamma_{k,0}$ intersects $C_k$ at $u'_{k,0}$. Here $j \in J_k$ takes integer values increasing with $t_{k,j}$, where $s_{k,j} = \sigma_{k,j} + it_{k,j}$ are non trivial zeros of $L(s\ ;\chi)$ belonging to $\Gamma_{k,j}$.

**Theorem 3**. *If $\chi$ is a primitive character, then there are infinitely many components $\Gamma'_m$ of the pre-image by $L(s\ ;\chi)$ or $\mathbb{R}$ which are mapped bijectively by $L(s;\ \chi)$ onto the interval $(1, +\infty)$ of the real axis. Consecutive curves $\Gamma'_m$ and $\Gamma'_{m+1}$ have no common points (except possibly for two of them) and they form strips $S_m$ which are mapped (not necessarily bijectively) onto the whole $z$-plane with a slit alongside the interval $(1, +\infty)$ of the real axis.*

*Proof*: Fig. 4 below illustrates the pre-image by $L(7,2,s)$ of an orthogonal mesh formed with rays starting at the origin and circles centered at the origin. Different colors have been used for the annuli so formed with saturation increasing counterclockwise and brightness increasing outward (the saturation is determined by the argument of the point and the brightness by its modulus) and the same color, brightness and saturation have been imposed to every point of the pre-image in order to facilitate the location of the respective points. The curves $\Gamma_{k,j}$ containing the non trivial zeros, as well as the curves $\Gamma'_k$ bounding the strips $S_k$ are all obvious. Some other curves crossing the real axis at branch points can be seen. Those branch points alternate with trivial zeros situated on the negative real half axis. Due to the fact that $L(s\,;\chi)$ is conformal in every domain which does not contain branch points, the quadrilaterals appearing in this figure, except those containing branch points, are conformally mapped by $L(s\,;\chi)$ onto the corresponding quadrilaterals in the $z$-plane.

Let us deal first with the exception case mentioned in the theorem. If we denote by $b_1, b_2, \ldots$ the branch points on the negative real half axis in the increasing order of their module, then at all the points $b_k$ we have configurations as those in [1], page 133, with $n = 2$. Two of the arcs starting at $b_1$ can be continued indefinitely along the interval $[L(b_1;\chi),+\infty)$ providing unbounded curves $\Gamma''_1$ in the upper half plane and $\Gamma''_{-1}$ in the lower half plane. The two curves are bijectively mapped by $L(s\,;\chi)$ onto the interval $[L(b_1;\chi),+\infty)$. Let us denote $\Gamma'_1 = \Gamma''_1 \cup [L(b_1;\chi),+\infty)$, respectively $\Gamma'_{-1} = \Gamma''_{-1} \cup [L(b_1;\chi),+\infty)$. Obviously, $\Gamma'_1 \cap \Gamma'_{-1} = [L(b_1;\chi),+\infty)$. Let us show that other consecutive curves $\Gamma'_k$ have no common point. Indeed, suppose that $s_0 \in \Gamma'_m \cap \Gamma'_n$. Then the sub-arcs of $\Gamma'_m$ and $\Gamma'_n$ corresponding to the interval of the real axis between $x = 1$ and $x = L(s_0;\,\chi)$ would bound a domain which is mapped by $L(s;\,\chi)$ onto the whole complex plane with a slit alongside this interval. Such a domain must contain a pole of $L(s;\,\chi)$, yet since $\chi$ is a primitive character, $L(s;\,\chi)$ is an integer function. Therefore, consecutive curves $\Gamma'_m$ and $\Gamma'_{m+1}$ , $m \in \mathbb{Z}\backslash\{0\}$ bound infinite strips $S_m$. Since the boundary of such a strip is mapped by $L(s;\,\chi)$ onto the interval $(1,+\infty)$, the strip itself is mapped, not necessarily bijectively, onto the whole complex plane with a slit alongside this interval. Suppose that $\Gamma'_{n_0}$ is the last curve of this type. Then the closed domain above it is mapped by $L(s;\,\chi)$ onto the whole complex plane. Indeed, let $s_0$ be a point in this domain and let $z_0 = L(s_0;\,\chi)$. Let us connect any point $z$ in the plane with $z_0$ by an arc $\gamma$ disjoint of the interval $(1,+\infty)$, except possibly for its end and which does not contain any branch point of $L(s\,;\chi)$. Continuation along $\gamma$ from $s_0$ brings us to a point $s$ such that $L(s;\chi) = z$. Thus, the whole $z$-plane is the image of that closed domain and $\Gamma'_{n_0}$ is mapped bijectively onto the interval $(1,+\infty)$. The neighborhood of a point on $(1,+\infty)$ must be the image of two half-neighborhoods of points on $\Gamma'_{n_0}$, which contradicts the injectivity of $L(s;\chi)$ on $\Gamma'_{n_0}$. A similar argument can be used if we suppose that there are no curves $\Gamma'_m$ below a $\Gamma'_{m_0}$. Therefore, we can let $m$ vary from $-\infty$ to $+\infty$ and count the consecutive curves $\Gamma'_m$ such that $S_0$ is the strip containing the real axis.

The way the strips $S_m$ have been formed shows that every one contains a unique unbounded component $C_m$ of the pre-image of the unit circle. Since every $C_m$ intersects the critical strip and $\Gamma'_m$ is situated between $C_m$ and $C_{m+1}$ we conclude that $\Gamma'_m$ intersects also the critical strip.

Obviously, all these components come in contact at $s = \infty$. If we increase the radius past 1, i.e. we take the pre-image of a circle $\gamma_r$ with $r > 1$, they all fuse into a unique unbounded curve intersecting every $\Gamma'_k$ as illustrated in Fig. 4 below.

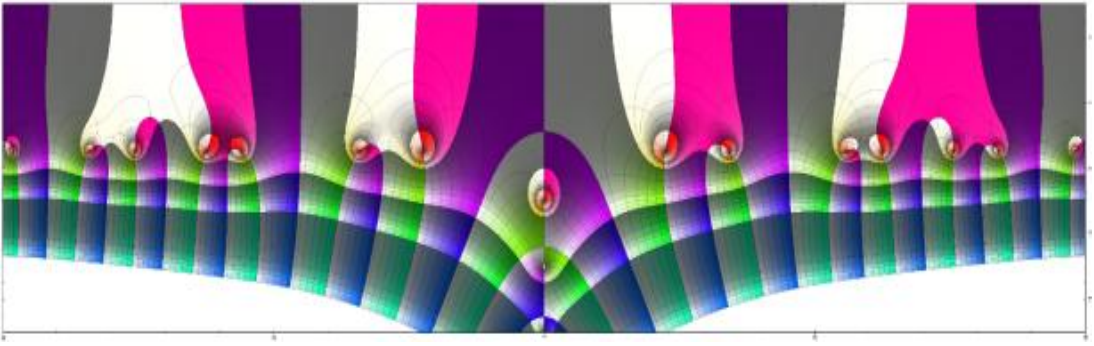

Fig. 4     Colored quadrilaterals conformally mapped by $L(7,4,s)$

Let us notice that if $\chi$ is real, then $\overline{L(s;\ \chi)} = L(\bar{s};\ \chi)$ and therefore the pre-image of the real axis is symmetric with respect to the real axis. Moreover, the real axis is mapped by $L(s;\ \chi)$ onto itself (see Fig. 4). Otherwise, there is no symmetry with respect to the real axis and this axis is not mapped onto itself (see Fig 5) below. However, in both cases, the trivial zeros are all real if $\chi$ is real and some are imaginary if $\chi$ is complex.

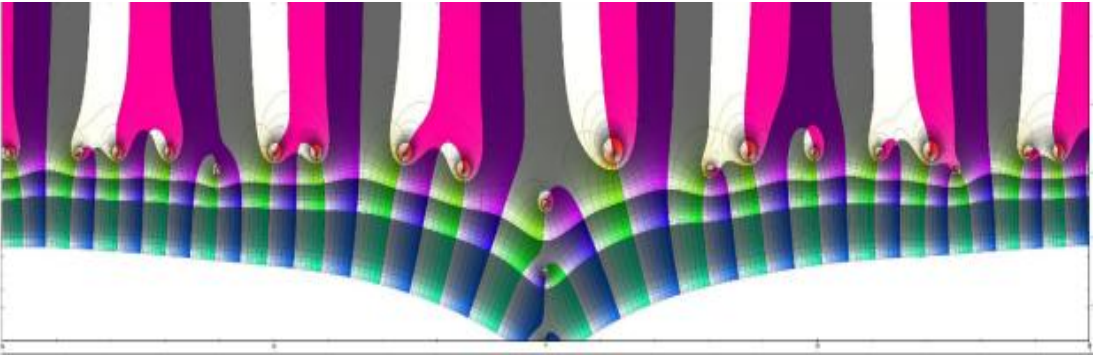

Fig. 5     Colored quadrilaterals conformally mapped by $L(7,2,,s)$

**Theorem 4**. *Every strip $S_k$ contains at least one non trivial zero of $L(s;\ \chi)$ and the number of non trivial zeros from $S_k$ is finite.*

Proof: Let us notice first that all the zeros of $L(s;;\ \chi)$ are interior to some strips since $0 \notin (1,+\infty)$ and therefore the pre-image of $0$ cannot belong to the pre-image of $(1,+\infty)$. Next, since $S_k$ is mapped onto the whole plane with the slit $(1,+\infty)$, the point zero must be the image of at least one point from $S_k$, which is a zero of $L(s;\chi)$. Finally, it is known that the non trivial zeros of $L(s;\ \chi)$ belong to the critical strip (see [14]) whose intersection with $S_k$ is a compact set. This set can contain only a finite number of zeros of the analytic function $L(s;\ \chi)$. However, the strip $S_0$ contains infinitely many trivial zeros of $L(s;\ \chi)$.

### 3. **Fundamental Domains of** $L(s; \chi)$

It is enough to deal with an arbitrary strip $S_k$. We will show later that all the zeros of $L(s; \chi)$ are simple, but for the moment, when referring to the number of zeros $s_{k,j}$ of $L(s; \chi)$ in $S_k$, we tacitly count zeros with their multiplicities. Suppose that their number is $j_k$. By considering them as leafs of a complete binary tree, whose internal nodes $v_{kj}$ are the zeros of the derivative of $L(s; \chi)$, and having in view the way these zeros were previously obtained, we conclude easily that the number of the branch points of $L(s; \chi)$ situated in $S_k$ is $j_k - 1$.

Let us notice that, since a point traveling in the same direction on a circle $\gamma_r$ intersects consecutively the positive and the negative real half axis, any one of the components of its pre-image should intersect consecutively the pre-image of the positive, respectively negative real half axis. If we color differently these last pre-images, then the colors such a point visits must alternate. We call this simple topological fact *the color alternating rule.* We met this phenomenon in Theorem 2 when dealing with the pre-image of the unit circle. The color alternating rule implied in that case that the points $u_{k,j}$ and $u'_{k,j}$ which were mapped by $L(s; \chi)$ respectively at $z = 1$ and $z = -1$ have been visited consecutively by a point $s$ turning around a component of the pre-image of the unit circle. When the radius of $\gamma_r$ increases past 1, all the unbounded components of this new circle cross the curves $\Gamma'_k$ (since the circle $\gamma_r$ intersects the interval $(1, +\infty)$ to whose pre-image every $\Gamma'_k$ belongs) and unite at the pre-image of $z = r$ into a unique unbounded component as shown in Fig. 4 and 5. Indeed the junction points are those points on $\Gamma'_k$ which are mapped by $L(s; \chi)$ at $z = r$. Continuing to increase $r$, the respective unbounded component will incorporate all the bounded components of the pre-image of $\gamma_r$ ( from any bounded region of the plane) after meeting them at the remaining branch points of $L(s; \chi)$ and the pre-image of $\gamma_r$ is a unique curve separating that region into two parts: the right one mapped by $L(s; \chi)$ inside the circle $\gamma_r$ and the left one mapped by $L(s; \chi)$ outside $\gamma_r$. The color alternating rule means that a point traversing this curve in the same direction will visit alternatively components of the pre-image of $\mathbb{R}$ having different colors. As this alternation happens for any $r$, and for $r$ big enough, there must be points $s = \sigma + it$ in $\Gamma_{k,j}$ with $s$ in the left region, we conclude that $\sigma \to -\infty$ when $s = \sigma + it \to \infty$, $s \in \Gamma_{k,j}$ whatever the color of $s$ would be and therefore the branches of $\Gamma_{k,j}, j \neq 0$ of different colors are adjacent. Moreover, since every $\Gamma'_k$ is situated between a $\Gamma_{k,j}$ and a $\Gamma_{k+1,j'}$ we have that $\sigma \to -\infty$ also when $s = \sigma + it \in \Gamma'_k$ and $L(s; \chi) \to +\infty$. It is obvious that two components of different colors can meet each other only at $\infty$ and at zeros and poles of $L(s; \chi)$ since $L(s; \chi)$ is an uniform function. Also it is impossible for two components of the same color to intersect each other, except if they start from a multiple zero or pole and we'll see in Section 4 that the functions $L(s; \chi)$ have no multiple zeros and the only pole can be $s = 1$, which when it exists, is simple. This analysis allows one to conclude that a similar landscaping to that from Fig. 4 and 5 should be found in any region of the plane and for any function $L(s; \chi)$.

**Theorem 5**. *If the strip $S_k$ contains $j_k$ non trivial zeros of $L(s; \chi)$, then it can be*

*divided into $j_k$ sub-strips which are fundamental domains of the function.*

Proof: Let us connect the images $L(v_{k,j}; \chi)$ of every branch point $v_{k,j}$ with $z = 1$ by a segment of line $\eta_{k,j}$ and perform continuations along every $\eta_{k,j}$ from the corresponding $v_{k,j}$. We obtain arcs and unbounded curves which, together with the components of the pre-image of the interval $(1, \infty)$, divide the strip $S_k$ into $j_k$ sub-strips. Indeed, as noticed in Section 2, $S_k$ contains $j_k - 1$ points $u_{k,j}$ such that $L(u_{k,j}; \chi) = 1$. The pre-image of $\eta_{k,j}$ are either arcs connecting two such points via $v_{k,j}$ or unbounded curves starting at $u_{k,j}$ and containing a point $v_{k,j}$. Together with the components of the pre-image of the interval $[1, +\infty)$ contained in the closure $\overline{S}_k$ of $S_k$ they bound infinite open strips $\Omega_{k,j}$ which are mapped conformally by $L(s; \chi)$ onto the complex plane with slits alongside the interval $[1, +\infty)$ plus $\eta_{k,j}$. Every component of the pre-image of the interval $[1, +\infty)$, except the boundary of $S_k$, belongs to two adjacent $\Omega_{k,j}$ and therefore the number of these fundamental domains (see [1], p.99) is $j_k$.

On the other hand, consecutive curves belonging to the pre-image of $\mathbb{R}$ and crossing the negative real half axis at $b_k$ (see Fig. 4) bound strips which are conformally mapped by $L(s; \chi)$ onto the complex plane with slits alongside $(-\infty, L(b_{2k})) \cup (L(b_{2k-1}), +\infty)$, $k = 1, 2, \ldots$ For a complex Dirichlet character $\chi$ (see Fig. 5) the domains between consecutive components of the pre-image of $\mathbb{R}$ containing trivial zeros, contain also branch points of $L(s; \chi)$, which can be located similarly to $v_{k,j}$. Then we follow the same procedure as in the case of the strips $S_k$ to find the corresponding fundamental domains.

### 4. The Derivative of $L(s; \chi)$

We have seen that, if $\chi$ is a Dirichlet non principal character, then $L(s; \chi)$ is an integer function, thus $L'(s; \chi)$ is also integer. The formula (2) implies that

(9)     $L'(s; \chi) = -\sum_{k=2}^{\infty} \frac{\chi(k) \ln k}{k^s}$

from which it is obvious that

(10)     $\lim_{\sigma \to +\infty} L'(\sigma + it; \chi) = 0$

Following the same reasoning as in the case of $L(s; \chi)$ we find that there are infinitely many components $\Upsilon'_k$ of the pre-image by $L'(s; \chi)$ of the real axis which are mapped bijectively by $L'(s; \chi)$ onto the interval $(-\infty, 0)$. They are mutually disjoint and consecutive components $\Upsilon'_k$ and $\Upsilon'_{k+1}$ border infinite strips $\Sigma_k$ which are mapped by $L'(s; \chi)$ not necessarily bijectively onto the complex plane with a slit alongside the interval $(-\infty, 0)$. Every $\Sigma_k$ contains a unique component $\Upsilon_{k,0}$ of the pre-image of $\mathbb{R}$ which is mapped bijectively onto the interval $(0, +\infty)$ and a finite number of components $\Upsilon_{k,j}, j \neq 0$ which are mapped bijectively onto $\mathbb{R}$. On each one of these last components $\operatorname{Re} s \to -\infty$ as $L'(s; \chi) \to \pm\infty$. The color alternation rule remains in force also for the pre-image by $L'(s; \chi)$ of $\mathbb{R}$.

**Theorem 6**. *The pre-image by $L'(s; \chi)$ of $\mathbb{R}$ intersects the pre-image by $L(s; \chi)$ of $\mathbb{R}$ at points where the tangent to the last one is horizontal and every point of this last pre-image where the tangent is horizontal is located at the intersection of the two pre-images.*

*Proof:* The mapping of every $\Omega_{k,j}$ can be obviously extended by continuity to the boundary of $\Omega_{k,j}$. Then, there is a natural parametrization $s = s(x)$ of every curve $\Gamma'_k$ and $\Gamma_{k,j}$ such that $L(s(x); \chi) = x$. Thus $L'(s(x); \chi) s'(x) = 1$, which implies that

(11)     $\arg L'(s(x); \chi) + \arg s'(x) = 0 \pmod{2\pi}$

In particular $\arg L'(s(x); \chi) = 0$ if and only if $\arg s'(x) = 0$ and $\arg L'(s(x); \chi) = \pi$ if and only if $\arg s'(x) = \pi$. The values $0$ and $\pi$ for $\arg L'(s(x_0); \chi)$ means that $s(x_0)$, which belongs to the pre-image by $L(s; \chi)$ of the real axis, is located on the pre-image by $L'(s; \chi)$ of the real axis, while the value $0$ and $\pi$ for $\arg s'(x_0)$ means that the tangent to the curve $x \to s(x)$ at $s(x_0)$ is horizontal. Vice-versa, if the tangent at $s(x_0)$ to the pre-image by $L(s; \chi)$ is horizontal, then $\arg s'(x_0)$ is $0$ or $\pi$, therefore $\arg L'(s(x_0); \chi)$ is $0$ or $\pi$, hence $s(x_0)$ belongs also to the pre-image by $L'(s; \chi)$ of the real axis. We call *intertwined* the corresponding components of the pre-image by the two functions of the real axis.

Let us notice that the projection $\gamma$ of the Riemann sphere of a curve $\Gamma'_k$ or $\Gamma_{k,j}$ is a smooth curve passing through the North Pole $N$ of the sphere which intersects the image $C$ on the sphere of the real axis only at $N$. The motion of the axis upward or downward corresponds to a contraction of $C$ keeping $N$ fixed and conserving the angle it makes with $\gamma$ at $N$. Since $\gamma$ is surrounded by $C$ in an obvious sense, there will be at least one point at which $\gamma$ and the contraction of $C$ touch each other. At the corresponding point in the complex plane the tangent to $\Gamma'_k$, respectively to $\Gamma_{k,j}$ is horizontal. Consequently, every one of these curves intertwines with a component of the pre-image of $\mathbb{R}$ by $L'(s; \chi)$. The following theorem shows how this happens.

**Theorem 7**. *Every curve $\Gamma'_k$ intertwines with a unique curve $\Upsilon'_k$ and with no curve $\Upsilon_{k,j}$ and every curve $\Upsilon'_k$ intertwines with a unique curve $\Gamma'_k$ and with no curve $\Gamma_{k,j}$.*

*Proof:* Let us denote by $a$ and $b$ the color of the components of the pre-image by $L(s; \chi)$ of the negative, respectively positive real half axis and by $c$ and $d$ the color of the components of the pre-image by $L'(s; \chi)$ of the negative, respectively positive real half axis. The curves $\Gamma'_k$ belong to the pre-image by $L(s; \chi)$ of the positive real half axis and therefore have color $b$, while $\Upsilon'_k$ belong to the pre-image by $L'(s; \chi)$ of the negative real half axis and therefore have color $c$. The theorem implies that, in what concerns the curves $\Gamma'_k$ and $\Upsilon'_k$, the color $b$ meets always the color $c$. Let us show that this property is more general, namely that if $s_0 \in \Gamma_{k,j} \cap \Upsilon_{k,j}$ and $\operatorname{Re} s_0 < \operatorname{Re} s_{k,j}(0)$, where $s = s_{k,j}(x)$, $x \in \mathbb{R}$ is the parametric equation of $\Gamma_{k,j}$ such that $L(s_{k,j}(x); \chi) = x$, hence $L'(s_{k,j}(x); \chi) s'_{k,j}(x) = 1$, then at $s_0$ the color $a$ meets the color $d$ or the color $b$ meets the color $c$. Indeed, suppose that $s = S_{k,j}(x)$, $x \in \mathbb{R}$ is the parametric equation of $\Upsilon_{k,j}$ such that $L'(S_{k,j}(x); \chi) = x$. There is $x_1, x_2 \in \mathbb{R}$ such that $s_0 = s_{k,j}(x_2) = S_{k,j}(x_1)$. Thus $L'(s_{k,j}(x_2); \chi) = L'(S_{k,j}(x_1); \chi) = x_1$, hence $x_1 s'_{k,j}(x_2) = 1$ therefore $\arg x_1 + \arg s'_{k,j}(x_2) = 0 \pmod{2\pi}$. We notice that $\operatorname{Re} s_0 < \operatorname{Re} s_{k,j}(0)$ implies that $\arg s'_{k,j}(x_2) = 0$ only if $x_2 < 0$ and $\arg s'_{k,j}(x_2) = \pi$ only if $x_2 > 0$. Then $x_1 > 0$ in the first case and $x_1 < 0$ in the second case, i.e. the color at $s_0$ of $\Gamma_{k,j}$ is $a$ iff that of $\Upsilon_{k,j}$ is $d$ and the color at $s_0$ of $\Gamma_{k,j}$ is $b$ iff that of $\Upsilon_{k,j}$ is $c$. We call these conditions *the color matching rule*. This color matching condition for $\Gamma'_k$ and $\Upsilon'_k$ is fulfilled with no restriction on $\operatorname{Re} s_0$, since their colors are respectively $b$ and $c$. We notice however that for $\operatorname{Re} s_0 > \operatorname{Re} s_{k,j}(0)$ this rule does not apply and its violation is obvious for $\Gamma_{k,0}$ which intersects in the part colored $b$ the curve $\Upsilon_{k,0}$ colored entirely $d$. We use the color matching rule in order to prove the theorem. The color alternating rule requires that $\Gamma_{k,j}$ adjacent to $\Gamma'_k$ or to $\Gamma'_{k+1}$ has the color $a$ on the side facing that component, therefore by the color matching rule the corresponding $\Upsilon_{k,j}$ should have the color $d$, which does not allow $\Upsilon_{k,j}$ to intersect $\Gamma'_k$ or $\Gamma'_{k+1}$ since in such a case the color matching rule would be violated: the color $d$ would meet the color $b$. Consequently, $\Gamma'_k$ can intertwine only with $\Upsilon'_k$ and vice-versa, and every $\Gamma_{k,j}$ can intertwine only with some $\Upsilon_{k,j'}$.

Fig. 6 and 7 below illustrate the correspondence between intertwined curves for $L(7, 4, s)$ and $L(7, 2, s)$.

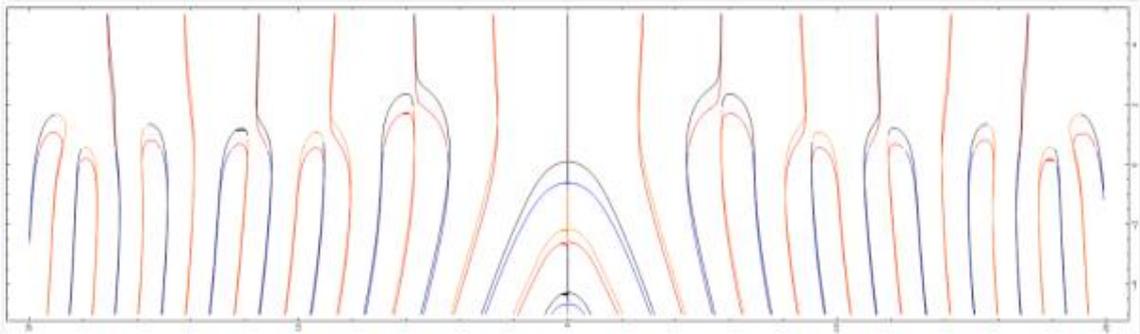

Fig.6     The correspondence between intertwined curves for $L(7, 4, s)$

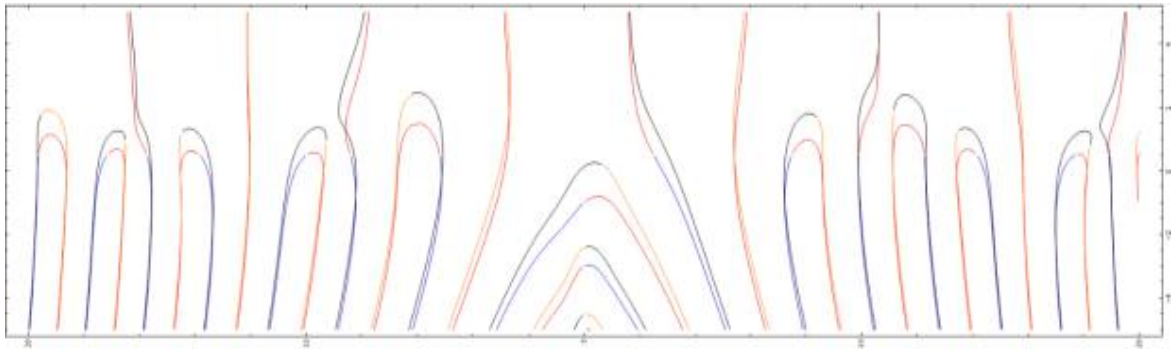

Fig.7     The correspondence between intertwined curves for $L(7, 2, s)$

### 5. The Multiplicity and the Location of the Zeros of $L(s \, ; \chi)$

It has been hypothesized (see [14]) that the zeros of every Dirichlet $L$-function are all simple. A proof of this affirmation can be based on the correspondence between intertwined curves. This correspondence and the conformal mapping properties of $L(s \, ; \chi)$ and of $L'(s \, ;)$ offer also the necessary tools for the proof of the Generalized Riemann Hypothesis.

**Theorem 8**. *The zeros of every Dirichlet $L$-series $L(s; \chi)$, as well as the zeros of its derivative are all simple.*

*Proof:* For the trivial zeros the affirmation is obvious. Suppose now that a non-trivial zero $s_{k,j}$ of an arbitrary $L(s; \chi)$ is multiple of order $m$. Then $m$ curves $\Gamma_{k,j}$ would intersect at $s_{k,j}$, which is a zero of order $m - 1$ for $L'(s; \chi)$. Thus, there are $m - 1$ curves $\Upsilon_{k,j}$ passing through $s_{k,j}$ and intertwining with $m - 1$ of the curves $\Gamma_{k,j}$. On the other hand, the $m$-th curve $\Gamma_{k,j}$ must have also an intertwining curve $\Upsilon_{k,j'}$. This last curve can meet any other $\Upsilon_{k,j}$ only at a zero, i.e. in this case at $s_{k,j}$. Then $s_{k,j}$ would be a multiple zero of order $m$ of $L'(s; \chi)$, which is a contradiction. A similar reasoning can be performed when studying the multiplicity of the zeros of $L'(s; \chi)$ and the conclusion of the theorem follows.

**Theorem 9**. *For every $q \geq 1$ and every Dirichlet character $\chi$ modulo $q$ the Riemann Hypothesis regarding the non trivial zeros of $L(s; \chi)$ is true: all the non trivial zeros have the real part equal to $1/2$.*

Proof: As already stated, we only need to deal with Dirichlet $L$-series $L(s; \chi)$ defined by primitive characters $\chi$. The Theorems 1-8 show that the pre-image of the real axis by such a function is similar to that produced by the Riemann Zeta function, except that when $\chi$ is not real, we have no symmetry with respect to the real axis. Moreover, the functional equation (7) allows us to use exactly the same arguments as those used for the Riemann Zeta function in order to generalize the Riemann Hypothesis. In the case of the Riemann Zeta function embraced curves $\Gamma_{k,j}$ and $\Gamma_{k,j'}$ have been found, as it appears in Fig. 8 below (see [13], p. 137 ).

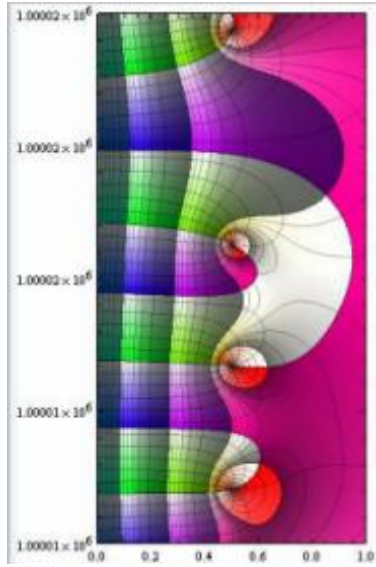

Fig. 8     Embraced curves $\Gamma_{k,j}$ exhibited by $\zeta$ in the range of $t = 10^6$

It is expected to find such configurations for some other Dirichlet $L$-functions. It can be easily checked that the zeros belonging to two such curves cannot have the same imaginary part. Indeed, let $I$ be the segment between the zeros belonging to the two curves. If $I$ is horizontal, then the color alternating rule requires that the arc between $u_{k,j}$ and $u_{k,j+1}$ which is projected  by $L(s\,;\chi)$ onto the corresponding cut in the $z$-plane intersects $I$. However, the images by $L'(s\,;\chi)$ of the two arcs must be disjoint, as shown in Fig. 9 below, or they should intersect of an even number of times, which is a contradiction.

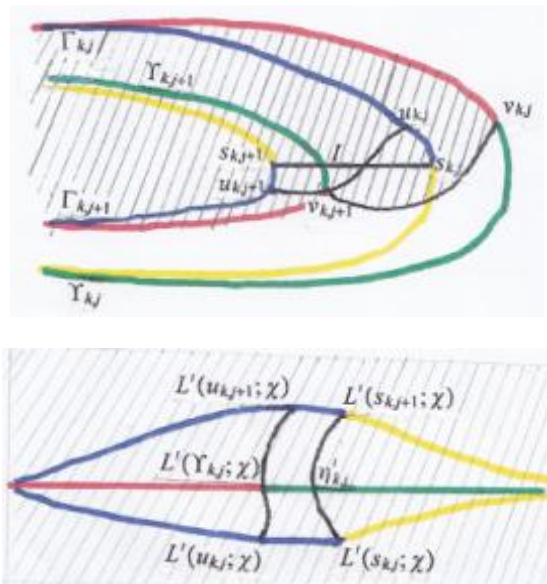

Fig.9 The zeros belonging to embraced curves cannot have the same imaginary part.

A similar contradiction appears if we take the image of $I$ by $L(s\,;\chi)$ when supposing that $\operatorname{Im} s_{k,0}$ coincide with $\operatorname{Im} s_{k,1}$ or with $\operatorname{Im} s_{k,-1}$. Therefore it is enough to deal only with zeros belonging to consecutive adjacent curves $\Gamma_{k,j-1}$ and $\Gamma_{k,j}$, $j \notin \{0,1\}$ included in the same strip $S_k$ or to $\Gamma_{k,j}$ and $\Gamma_{k+1,j'}$ included in consecutive strips $S_k$ and $S_{k+1}$. In each one of these two cases there are two distinct relative positions of the respective curves allowing $I$ to be horizontal. By inspecting the images of $I$ by $L(s\,;\chi)$ and by $L'(s\,;\chi)$ we are brought to contradictions in all these cases.

*Case (a).* Suppose that the respective zeros are $s_{k,j-1} = \sigma_{k,j-1} + it$ and $s_{k,j} = \sigma_{k,j} + it$. Let $\Gamma_{k,j-1}$ and $\Upsilon_{k,j-1}$, respectively $\Gamma_{k,j}$ and $\Upsilon_{k,j}$ be the corresponding intertwined curves with $s_{k,j} \in \Gamma_{k,j}$, $s_{k,j-1} \in \Gamma_{k,j-1}$, $v_{k,j} \in \Upsilon_{k,j}$ and $v_{k,j-1} \in \Upsilon_{k,j-1}$. The parametric equation of the segment $I$ is

$$(12) \qquad s(\lambda) = (1-\lambda)s_{k,j} + \lambda s_{k,j-1}, \qquad 0 \leq \lambda \leq 1.$$

The image of $I$ by $L(s\,;\chi)$ is a closed curve $\eta_{k,j}$ of equation

$$(13) \qquad z(\lambda) = L((1-\lambda)s_{k,j} + \lambda s_{k,j-1}\,;\chi), \qquad 0 \leq \lambda \leq 1$$

with $z(0) = z(1) = 0$. It is a smooth curve, except for the point $v_{k,j}$ if it happens that $v_{k,j} \in \eta_{k,j}$. Then we can differentiate with respect to $\lambda$ for $\lambda \in (0,1)\backslash\{v_{k,j}\}$ and we get

$$(14) \qquad z'(\lambda) = L'((1-\lambda)s_{k,j} + \lambda s_{k,j-1})(s_{k,j-1} - s_{k,j}),$$

which is an arc $\eta'_{k,j}$. Since $s_{k,j-1} - s_{k,j} = \sigma_{k,j-1} - \sigma_{k,j} > 0$, we have that

$$(15) \qquad \arg z'(\lambda) = \arg L'((1-\lambda)s_{k,j} + \lambda s_{k,j-1}),$$

which means that the angle between the tangent to $\eta_{k,j}$ at a point $z(\lambda)$, $\lambda \neq v_{k,j}$ and the positive real half axis coincide with that made by the position vector of $L'((1-\lambda)s_{k,j} + \lambda s_{k,j-1})$ and the positive real half axis. Since $\lim_{\lambda \searrow 0} z'(\lambda)$ and $\lim_{\lambda \nearrow 1} z'(\lambda)$ exist, the equality (15) holds also for the ends of $\eta_{k,j}$ and $\eta'_{k,j}$. The hypothetical configuration shown in Fig. 10.a below implies that the two half-tangents to $\eta_{k,j}$ at $z = 0$ point one to the upper half plane and the other to the lower half plane.

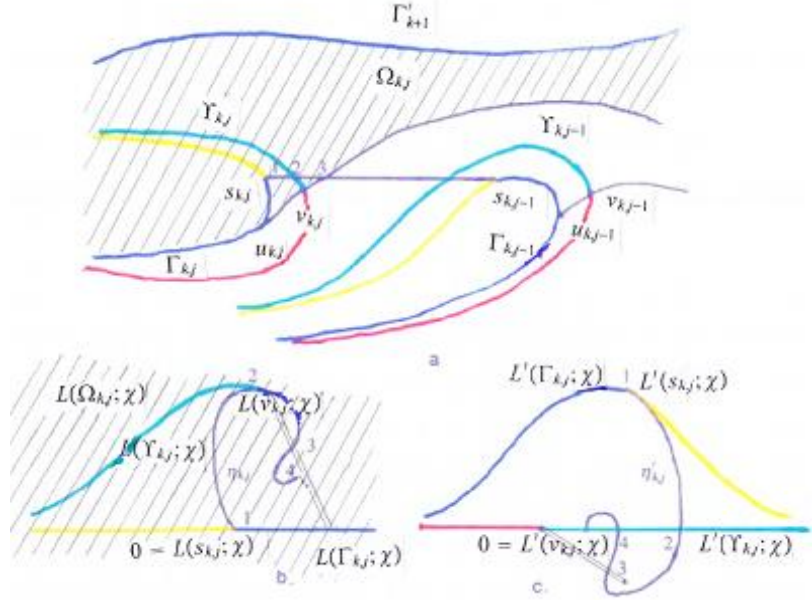

Fig. 10    Hypothetical situation with $\operatorname{Im} s_{k,j} = \operatorname{Im} s_{k,j-1}$

However, the arguments of $L'((1-\lambda)s_{k,j} + \lambda s_{k,j-1})$ at the ends of $\eta'_{k,j}$ are both between

$0$ and $\pi$, which contradicts the equality (15). Therefore a configuration like that in Fig. 12.a cannot exist.

*Case (b).* Let us deal with a hypothetical configuration like that shown in Fig. 11.a. This time the parametric equation of the segment $I$ is:

$$(16) \qquad s(\lambda) = (1-\lambda)s_{k,j} + \lambda s_{k,j+1}, \qquad 0 \le \lambda \le 1$$

and as in the previous case, the images of $I$ by $L(s;\chi)$ and by $L'(s;\chi)$ are respectively a closed curve $\eta_{k,j} : z = z(\lambda),\ 0 \le \lambda \le 1$ passing through the origin and having well defined half-tangents at the origin and an arc $\eta'_{k,j}$ starting at $L'(s_{k,j};\chi)$ and ending at $L'(s_{k,j+1};\chi)$. Since

$$(17) \qquad \arg z'(\lambda) = \arg L'(s(\lambda);\chi),$$

when $\arg z'(\lambda) = 0$, the arc $\eta'_{k,j}$ intersects the positive real half axis and when $\arg z'(\lambda) = \pi$ the arc $\eta'_{k,j}$ intersects the negative real half axis. Let us denote by $1, 2, 3, 4, \ldots$ some corresponding points on $I$ and respectively on $\eta_{k,j}$ and $\eta'_{k,j}$ as shown in Fig.11.

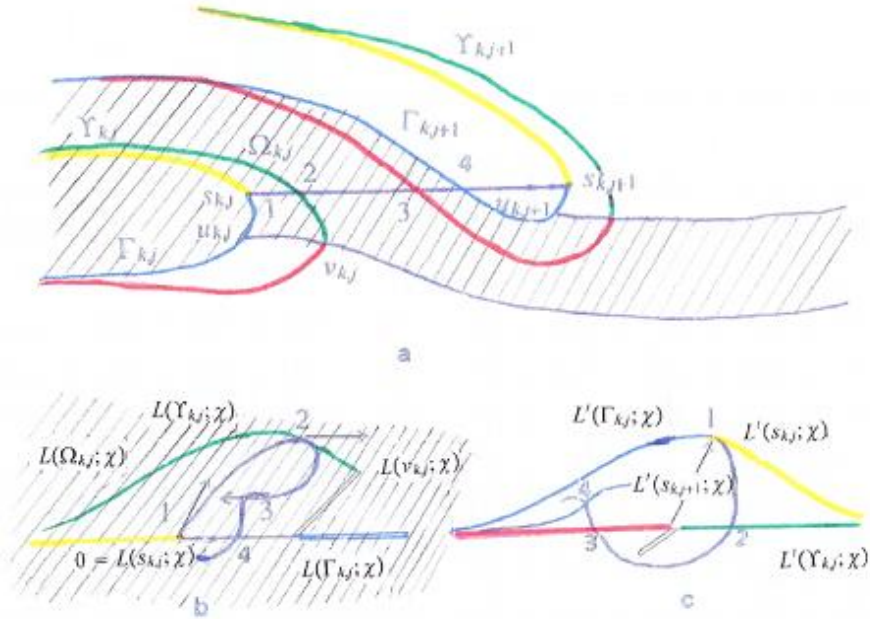

Fig. 11    Hypothetical situation in which $\operatorname{Im} s_{k,j} = \operatorname{Im} s_{k,j+1}$

We notice that after $2$ the curve $\eta_{k,j}$ must turn to the origin and therefore $\arg z'(\lambda) = \pi$ somewhere at a point $3$ and and $\eta'_{k,j}$ must intersect the negative real half axis at the corresponding point $3$. If $\eta_{k,j}$ reaches directly the origin, without intersecting the positive real half axis at the point $4$ as seen in Fig. 11b, then $I$ must reach $s_{k,j+1}$ without intersecting $\Gamma_{k,j+1}$ at $4$ as shown in Fig. 11a, i.e. $4$ is in fact $s_{k,j+1}$. Then $\eta'_{k,j}$ must end in $L'(s_{k,j+1}; \chi)$ which belongs to the upper half plane, while the corresponding half-tangent at $0$ to $\eta_{k,j}$ points to the lower half plane, which contradicts (17). If $\eta_{k,j}$ does not reach the origin, but intersects the positive real half axis in $4$, as seen in Fig.11b, then at a point $5$ we should have $\arg z'(\lambda) = \pi$, which forces $\eta'_{k,j}$ to intersect again the negative real half axis. Then, when $\eta_{k,j}$ reaches the origin, $\eta'_{k,j}$ would be still in  the lower half plane, which contradicts again the equality (17). In other words a configuration like that shown  in Fig.11 is impossible.

Case (c). Suppose now that the consecutive zeros $s_{k,j'}$ and $s_{k+1,j}$ belonging respectively to the strips $S_k$ and $S_{k+1}$ are of the form $s_{k+1,j} = \sigma + it$ and $s_{k,j'} = 1 - \sigma + it$, $0 < \sigma < 1/2$, as shown in Fig. 12a below.

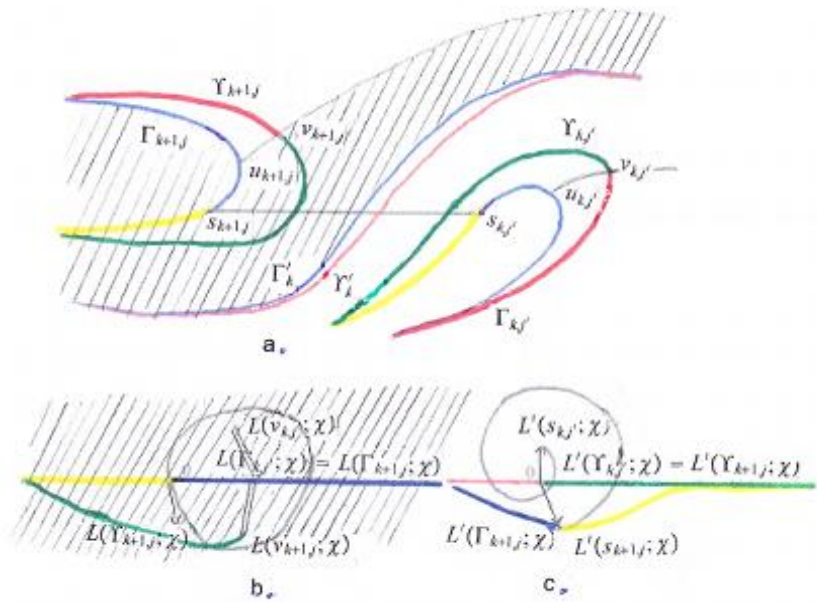

Fig . 12    Hypothetical situation in which $\operatorname{Im} s_{k+1,j} = \operatorname{Im} s_{k,j'}$

We notice that $L'(s_{k,j'}; \chi)$ belongs to the upper half plane, while $L'(s_{k+1,j}; \chi)$ to the lower half plane. The half tangents at $0$ to the image by $L(s ; \chi)$ of the interval $I$ between the two zeros point both to the lower half plane and having in view the similarities with the previous cases, we conclude that a configuration as that in  Fig. 13 is impossible.

*Case (d).* Suppose now that $s_{k,j} = \sigma + it$ and $s_{k+1,j'} = 1 - \sigma + it$, $0 < \sigma < 1/2$, as shown in Fig. 13a below.

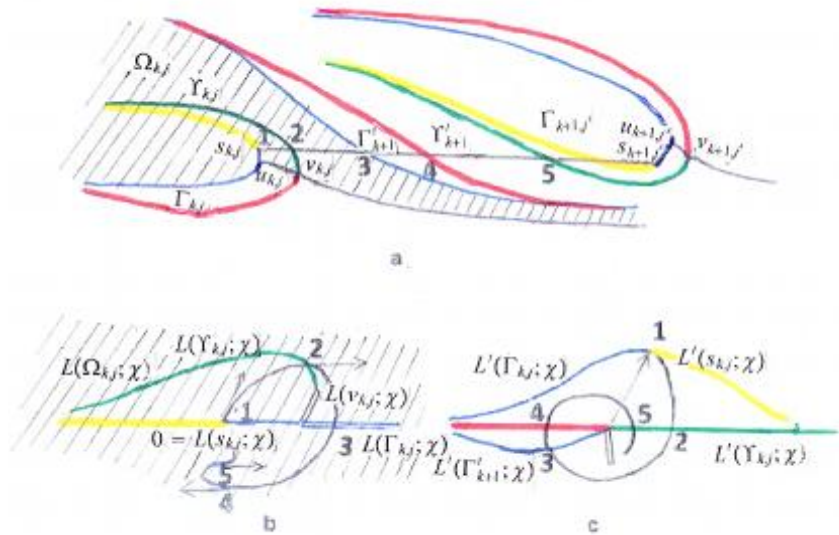

Fig. 13    Hypothetical situation in which $\operatorname{Im} s_{k,j} = \operatorname{Im} s_{k+1,j'}$

Using the same notations as in the previous cases, we notice that after the position 5 the vector $L'(s(\lambda) ; \chi)$ points to the lower half plane. If $\eta_{k,j}$ turns after position 5 directly to the origin, then $z'(\lambda)$ will point to the upper half plane, which contradicts the formula (17). If it crosses the negative half axis and then turns back to the origin, then $\arg z'(\lambda) = 0$ for some value $\lambda$ and $\eta'_{k,j}$ will cross the negative half axis ending in the upper half plane, which contradicts again (17). Therefore, a configuration like that shown in Fig 13a is impossible.

We conclude that, due to the functional equation (7), all the non trivial zeros of $L(s ; \chi)$ in $S_k$ have the real part $1/2$, for any $k$, i.e. the Riemann Hypothesis generalizes to the class of Dirichlet $L$-functions.

**Acknowledgments**: The author is grateful to Florin Alan Muscutar for fruitful discussions on the topics of this paper and for providing computer generated graphics.

**References**

[1] Ahlfors, L. V., *Complex Analysis*, International Series in Pure and Applied Mathematics, 1979

[2] Ahlfors, L. V. and Sario, L., *Riemann Surfaces,* Princeton University Press, 1960

[3] Andreian Cazacu, C. and Ghisa, D., *Global Mapping Properties of Analytic Functions*, Proceedings of 7-th ISAAC Congress, London, U.K., 3-12, 2009

[4] Andreian Cazacu, C. and Ghisa, D., *Fundamental Domains of Gamma and Zeta Functions,* IJMMS, ID 985323, 21 pages, 2011

[5] Ballantine, C and Ghisa, D., *Color Visualization of Blaschke Product Mappings*, Complex Variables and Elliptic Equations, Vol. 55, Nos. 1-3, 201-217, 2009

[6] Ballantine, C. and Ghisa, D., *Global Mapping Properties of Rational Functions,* Proceedings of 7-th ISAAC Congress, London, U.K., 13-22, 2009

[7] Barza, I. and Ghisa, D. , *The Geometry of Blaschke Product Mappings,* H.G.W.Begehr, A.O.Celebi and R.P.Gilbert eds., Further Progress in Analysis, World Scientific, 197-207, 2008

[8] Barza, I. and Ghisa, D., *Blaschke Product Generated Covering Surfaces*, Mathematica Bohemica, 134, No. 2, 173-182, 2009

[9] Cao-Huu, T. and Ghisa, D., *Invariants of Infinite Blaschke Products*, Mathematica, Tome 49 (720) No. 1, 24-34, 2007

[10] Ghisa, D., *The Riemann Hypothesis,* Studia Univ. Babes-Bolyai, 57, No. 2, 195-208, 2012

[11] Ghisa, D. and Muscutar, F. A., *The Riemann Zeta Function Has Only Simple Zeros,* International J. of Math. Sci. & Engg. Appls., Vol. 6, No III 239-250, 2012

[12] Ghisa, D.  and Salagean,  G. S., *Distortion and Related Inequalities for the Riemann Zeta Function*, International J. of Sci. & Engg. Appls. Vol. 6, No. IV, 303-315, 2012

[13] Ghisa, D., *Fundamental Domains and the Riemann Hypothesis,* Lambert Academic Publishing, 2013

[14] Montgomery, H. L. and Vaughan, R. C.,  *Multiplicative Number Theory,* Cambridge Studies in Advanced Mathematics 97, 2007